\documentclass{article}%
\usepackage{amsmath}%
\usepackage{amsfonts}%
\usepackage{amssymb}%
\usepackage{graphicx}%
\usepackage{color} 

\newtheorem{thm}{Theorem}

\newtheorem{cor}[thm]{Corollary}

\def\H{{\cal H}}

\def\Bbb{\mathbb}

\def\l{\lambda}

\def\per{\mathop{\rm per}}

\def\proof{{\noindent \bf Proof. \hspace{.01in}}}
\newcommand{\qed}{\hspace{.1in} \vrule height 7pt width 5pt depth 0pt \medskip}

\begin{document}
\title{\Large An inequality for tensor product of positive
operators and its applications}
\author{Haixia Chang$^{\,\rm a,}$\thanks{{\em E-mail address:} hcychang@163.com.},\;
Vehbi E. Paksoy$^{\,\rm b,}$\thanks{{\em E-mail address:} vp80@nova.edu.}, \;
Fuzhen Zhang$^{\,\rm b,}$\thanks{{\em E-mail address:} zhang@nova.edu.}\\
\footnotesize{$\mbox{}^{\rm a}$ Department of Applied Mathematics, Shanghai Finance University}\\
\footnotesize{Shanghai 201209, P.R.China}\\
\footnotesize{$\mbox{}^{\rm b}$ Farquhar College of Arts and Sciences,
Nova Southeastern University}\\
\footnotesize{3301 College Ave., Fort Lauderdale, FL 33314, USA}
}
\date{}
\maketitle

\noindent In memory of  Hans Schneider

\bigskip
 \hrule
\bigskip

\noindent {\bf Abstract}

\medskip
We present an inequality for tensor product of positive operators on Hilbert spaces by considering the tensor product of operators as words on certain alphabets (i.e., a set of letters).  As applications of the operator inequality and by a multilinear approach,
we show some matrix inequalities concerning induced operators and generalized matrix functions
(including determinants and permanents as special cases).

\medskip
{\footnotesize
\noindent {\em AMS Classification:} 47A63, 15A48

\noindent {\em Keywords:} generalized matrix function, induced operator,
inequality,
positive operator, positive semidefinite matrix,  positivity, tensor, word}

\bigskip
\hrule


\section{Introduction}

Let $\H$ be a Hilbert space over the complex number field $\Bbb C$ with an inner product $\langle \cdot, \cdot\rangle $.  Denote by ${\cal B}(\H)$ the ${\rm C}^*$-algebra of all bounded linear operators on $\H$.
We write $A\geq 0$ if $A$ is a  positive semidefinite operator  on $\H$ (we simply call it a
positive operator), that is, $A$ is self-adjoint and $\langle Ax, x\rangle \geq 0$ for all $x\in \H$.
For self-adjoint $A, B\in {\cal B}(\H)$, we write $A\geq B$ if $A-B\geq 0$.
It is well known that if $A\geq 0$ and $B\geq 0$ then the sum $A+B\geq 0$ (on $\H$) and tensor product $A\otimes B\geq 0$ (on $\otimes^2 \H=\H \otimes \H$).
Moreover, if $A$ is positive then
 the tensor product $\otimes^m A =A\otimes \cdots \otimes A $ ($m$ copies of $A$) is positive (on $\otimes ^m \H$)  for any positive integer $m$.
For the finite-dimensional case of $\H$, we denote by $\l_{\min}(A)$ and $\l_{\max}(A)$ the smallest and largest eigenvalues of the positive linear
operator (matrix) $A$ on $\H$, respectively.


In Section 2, we present an inequality for
positive operators with tensor product. In Section 3, as applications of our main result, we  deduce inequalities for generalized matrix functions, including determinants and permanents. Our results may be regarded as additions to the
recent ones in the research development of positivity (see, e.g., \cite{Bhatia14, Gao14, StanfordAIM14, AMSMeeting15, DBalaAMSTrans14, LawL14, LiZ14, ZFZActa12}).

\section{Main results}

We present our main result in this section. The proof is accomplished  by using the idea of words
(which, for instance, is used to show normality of matrices in \cite{WZ95words}). Let $\{A_i, \dots, A_j\}$ be a
(multi-)set of operators from ${\cal B}(\H)$ in which $t$ operators are distinct. For example, $t=2$ for $\{A_1, A_1, A_2\}$.
A {\em tensor word} or  {\em word}  of operators $A_i, \dots, A_j$ on $\H$ of {length} $m$ and of $t$ (distinct) representatives with respect to the
 tensor product $\otimes$, symbolized by $w_m^{t}(A_i, \dots, A_j)$, or $w^t(A_i, \dots, A_j)$, or even simply
$w^t$ (if no confusion is caused), is a tensor product
$$A_{s_1}\otimes \cdots \otimes A_{s_m},$$
in which $A_{s_1}, \dots, A_{s_m}$ are taken from $\{A_i, \dots, A_j\}$, and among $A_{s_1}, \dots, A_{s_m}$,
there are $t$ distinct operators. For instance,
$A_1\otimes A_1=\otimes^2A_1$ is a $w^1$ word;  $A_1\otimes A_2$ is a
$w^2$ word, and $A_3\otimes A_2$ is also a $w^2$ word. Note that $w^2_3(A_1, A_2)$ may represent any of
$A_1\otimes A_1\otimes A_2$, $A_1\otimes A_2\otimes A_1$, $A_2\otimes A_1\otimes A_1$,
$A_1\otimes A_2\otimes A_2$, $A_2\otimes A_1\otimes A_2$, and $A_2\otimes A_2\otimes A_1$.
So when we say a  $w^2(A_1, A_2)$ word, we mean one of those tensor words (with the given length $m=3$).

\begin{thm}\label{Thm2}
Let $A_{1}, A_{2}, \dots,A_{k}\in {\cal B}(\H)$ be  positive operators. Then
\begin{eqnarray}
 \lefteqn{\otimes^{m}(A_{1}+A_{2}+\cdots+A_{k})} \label{Thm2eq1}\\
 & & -\sum_{\footnotesize {1\leq i_{1}<\dots<i_{k-1}\leq k}  } \otimes^{m}(A_{i_{1}}+A_{i_{2}}+\cdots+A_{i_{k-1}})\label{Thm2eq2}\\
 & & +\sum_{\footnotesize {1\leq j_{1}<\dots<j_{k-2}\leq k}  }  \otimes^{m}(A_{j_{1}}+A_{j_{2}}+\cdots+A_{j_{k-2}})-\cdots\label{Thm2eq3}\nonumber\\
  & & +(-1)^{k-t}\sum_{\footnotesize {1\leq x_{1}<\dots<x_{t} \leq k} }  \otimes^{m}(A_{x_{1}}+A_{x_{2}}+\cdots+A_{x_{t}})+\cdots\label{Thm2eq3b}\\
 & & +(-1)^{k-2}\sum_{\footnotesize {1\leq p_{1}<p_{2}\leq k }}  \otimes^{m}(A_{p_{1}}+A_{p_{2}})\label{Thm2eq4}\\
 & & +(-1)^{k-1}\sum_{\footnotesize {\substack{q=1}}}^k \otimes^{m} A_{q} \geq 0.\label{Thm2eq5}
\end{eqnarray}
\end{thm}

\proof The idea of the proof is to show that all the tensor words $w^t$ with  $t<k$ representatives
 are cancelled out in the reduced form (after the additions and subtractions). That is,
after the calculations, the only tensors that survive will be the ones of the form in the summation
$$\sum_{\footnotesize {\substack{ \mbox{ each of $A_1, A_2,\dots, A_k$ } \\ \mbox{  appears at least once }}}}
 ( A_{i_{1}}\otimes \cdots \otimes A_{i_{m}}).$$

 Consider  the words $w^1$, i.e., the tensors of one operator, $\otimes ^m A_i$, say $\otimes ^m A_1$.
 The tensor $\otimes ^m A_1$ appears in (\ref{Thm2eq1}) once, in (\ref{Thm2eq2}) $(k-1)={k-1\choose k-2}$ times, $\dots$,
 in (\ref{Thm2eq4}) $(k-1)={k-1\choose 1}$ times, and in (\ref{Thm2eq5}) once. Therefore, after cancelation,
 the  total number of $\otimes ^m A_1$ left in the reduced form is
 $$1- {k-1 \choose k-2}+\cdots +(-1)^{k-2}{k-1 \choose 1}+(-1)^{k-1}=(1-1)^{k-1}=0.$$

Now we consider a general tensor word $A_{i_{1}}\otimes \cdots \otimes A_{i_{m}}$. Let it be a $w^t$ word
with $t<k$, and let $A_{x_1}, A_{x_2}, \dots, A_{x_t}$ be  the distinct representatives of
 $A_{i_{1}}, A_{i_2}, \dots,  A_{i_{m}}$. We show that $A_{i_{1}}\otimes \cdots \otimes A_{i_{m}}$
  vanishes in the reduced form.
 Since it is a $w^t$ word, it cannot be a $w^{t-1}$ word. Therefore, this $w^t$ word does not appear in
  the expansion of any sum-tensor $\otimes^{m}(A_{y_{1}}+A_{y_{2}}+\cdots+A_{y_{s}})$, $s<t$. That is,
$\otimes^m (A_{x_1}+A_{x_2}+\cdots +A_{x_t})$ in (\ref{Thm2eq3b}) is the only sum-tensor containing
all $A_{x_1}, A_{x_2}, \dots, A_{x_t}$. It is important to observe that a $w^l$ word contains
all $A_{x_1}, A_{x_2}, \dots, A_{x_t}$ if and only if $l\geq t$ and $w^l$ is obtained from a
sum-tensor $\otimes^m (A_{x_1}+A_{x_2}+\cdots +A_{x_t}+A_u+\cdots +A_v)$, in which, after expansion,
each word $w^t(A_{x_1}, A_{x_2}, \dots, A_{x_t})$ appears once and only once; so does
$A_{i_{1}}\otimes \cdots \otimes A_{i_{m}}$. 

Note that for each fixed word $w^t(A_{x_1}, A_{x_2}, \dots, A_{x_t})$ (say, $A_{i_{1}}\otimes \cdots \otimes A_{i_{m}}$)
 in the sum-tensor
$\otimes^m(A_{x_1}+A_{x_2}+\cdots +A_{x_t})$, there are $k-t={k-t\choose 1}$ many terms
$\otimes^m (A_{x_1}+A_{x_2}+\cdots +A_{x_t}+A_x)$  in the expression prior to (\ref{Thm2eq3b});  they are
$$\otimes^m (A_{x_1}+A_{x_2}+\cdots +A_{x_t}+A_x), \quad x\in \{1, 2, \dots, k\}\setminus \{x_1, \dots, x_t\}.$$

Similarly,  considering one level above,
there are ${k-t\choose 2}$ many terms in
$$\sum \otimes^m (A_{x_1}+A_{x_2}+\cdots +A_{x_t}+A_x+A_y), \; x<y, \; x, y\in \{1, 2, \dots, k\}\setminus \{x_1, \dots, x_t\}$$ containing $\otimes^m (A_{x_1}+A_{x_2}+\cdots +A_{x_t})$.
Thus,  going in this way
up to
(\ref{Thm2eq1}), we have the total number of the $w^t(A_{x_1}, A_{x_2}, \dots, A_{x_t})$ word $A_{i_{1}}\otimes \cdots \otimes A_{i_{m}}$
 appearing in the
tensors $\otimes^m (A_{x_1}+A_{x_2}+\cdots +A_{x_t}+\cdots )$ in the reduced form:
$$1- {k-t \choose k-t-1}+\cdots +(-1)^{k-t-1}{k-t \choose 1}+(-1)^{k-t}=(1-1)^{k-t}=0.$$

It follows that the only remaining tensors  are  $w^k$ words, i.e., $A_{i_{1}}\otimes \cdots \otimes A_{i_{m}}$ in which every $A_1, A_2, \dots, A_k$
has to appear at least once (so $m\geq k$). That is, in the reduced form, namely, after computing  all expressions (\ref{Thm2eq1}) through (\ref{Thm2eq5}), we obtain
 $$\sum_{\footnotesize {\substack{ \mbox{ each of $A_1, A_2,\dots, A_k$ } \\ \mbox{  appears at least once }}}}
 ( A_{i_{1}}\otimes \cdots \otimes A_{i_{m}})  \geq 0.$$

If $m<k$,  then the left hand side of the above display is 0. \qed

The inequality  below is a special case of Theorem \ref{Thm2}. We single it out as it is
 interesting in its own right and all our results in Section 3 rely on it.

\begin{thm}\label{Thm1}
Let $A_{1}, A_{2}, A_{3}\in {\cal B}(\H)$ be positive operators.
  Then, for any positive integer $m$,
\[ \big (\otimes^{m}(A_1+A_2+A_3)+\otimes^{m}A_{1}+\otimes^{m}A_{2}+\otimes^{m}A_{3}\big )\]
\[  {}\quad -\big (\otimes^{m}(A_1+A_2)+\otimes^{m}(A_1+A_3)+\otimes^{m}(A_2+A_3)\big ) \geq 0.\]
If $\H$ is finite-dimensional,  then the eigenvalues   of the above difference lie between
$$\big ( 3(3^{m-1}-2^{m}+1) \big )\cdot \big ( \min_{1\leq i\leq 3} \l_{\min}(A_i) \big )^m$$
and
$$\big ( 3(3^{m-1}-2^{m}+1) \big )\cdot \big ( \max_{1\leq i\leq 3} \l_{\max}(A_i) \big )^m.$$
\end{thm}

\proof
 If $m=1$ or $m=2$, it is straightforward to check that the difference is zero.
  Let $m\geq 3$. From Theorem \ref{Thm2}, we have
 \begin{eqnarray*}
  \lefteqn{\otimes^{m}(A_1+A_2+A_3)+\otimes^{m}A_{1}+\otimes^{m}A_{2}+\otimes^{m}A_{3}}\\
  & & -\otimes^{m}(A_1+A_2)-\otimes^{m}(A_1+A_3) -\otimes^{m}(A_2+A_3)\\
&= &  \sum_{\footnotesize {\substack{ \mbox{ each of 1,2,3 } \\ \mbox{  appears at least once }}}}( A_{i_{1}}\otimes \cdots \otimes A_{i_{m}}) \quad \geq 0.
\end{eqnarray*}

 When the indices are taken from the set $\{1,2,3 \} $ and not all the same, there are $3^{m}-3$ different terms. When the indices are taken from the 2-element sets,
  each summation will consist of $2^{m}-2$ distinct tensor. Clearly, the first sum of tensors with indices  from $\{1,2,3\}$ contains all possible configurations of tensors in the other three sums. On the other hand, the first summation has additional tensors which contains all of the three elements from $\{1,2,3 \}$. So, the last sum in above computation has exactly $3^{m}-3 -3(2^{m}-2)=3(3^{m-1}-2^{m}+1)$ distinct positive semidefinite entry.

 For the smallest eigenvalue, since for every tensor product $A_{i_{1}}\otimes \cdots \otimes A_{i_{m}}$,
 $$\l_{\min}(A_{i_{1}}\otimes \cdots \otimes A_{i_{m}})=\prod_{t=1}^m \l_{\min}(A_{{i_t}}),$$
 we see that the smallest eigenvalue of the difference has a lower bound
 $$\big ( 3(3^{m-1}-2^{m}+1) \big )\cdot \big ( \min_{1\leq i\leq 3} \l_{\min}(A_i) \big )^m.$$
Similarly, we arrive at an upper bound for the largest eigenvalue:
 $$\big ( 3(3^{m-1}-2^{m}+1) \big )\cdot \big ( \max_{1\leq i\leq 3} \l_{\max}(A_i) \big )^m.
 \quad  \qed$$

\section{Applications}
Let $\H$ be finite-dimensional and $\otimes^m \H$ be the tensor product space of $m$ copies of $\H$.
Let $G$ be a subgroup of the $m$-symmetric group $S_m$ on $m$ letters and let $\chi $ be an irreducible character on $G$. Denote by
 $V_{\chi}(G)$ the symmetry class of tensors associated with $G$ and $\chi$
 (see, e.g., \cite[p.~154]{MerMul97}).  For a linear operator $A$ on $\H$, the induced operator
 $K(A)$ of $A$ with respect to $G$ and $\chi$ is the restriction of $\otimes^m A$ on $V_{\chi}(G)$; that is,
 $K(A)=(\otimes^m A)|_{V_{\chi}(G)}$ (see, e.g., \cite[p.~185, p.~235]{MerMul97}).

So, for the positive operators $A_1, A_2, A_3$ in Theorem \ref{Thm1} with $\dim \H<\infty$,
\begin{eqnarray*}\label{K}
(\otimes^{m}(A_1+A_2+A_3))|_{V_{\chi}(G)}+(\otimes^{m}A_{1})|_{V_{\chi}(G)}
+(\otimes^{m}A_{2})|_{V_{\chi}(G)}+(\otimes^{m}A_{3})|_{V_{\chi}(G)}\end{eqnarray*}
\vspace{-.25in}
\begin{eqnarray*}\label{K}
-\big (\otimes^{m}(A_1+A_2))|_{V_{\chi}(G)}+(\otimes^{m}(A_1+A_3))|_{V_{\chi}(G)}+(\otimes^{m}(A_2+A_3))|_{V_{\chi}(G)} \big )\geq 0,\end{eqnarray*}
 that is,
\begin{eqnarray*}\label{K}
 K(A_1+A_2+A_3)+K(A_{1})+K(A_{2})+K(A_{3})
 \end{eqnarray*}
\vspace{-.25in}
\begin{eqnarray*}\label{K}
\qquad- \big (K(A_1+A_2)+K(A_1+A_3)+K(A_{2}+A_{3})\big ) \geq 0.
 \end{eqnarray*}

The induced operators  are closely related to generalized matrix functions (see, e.g., \cite[p.~213]{MerMul97}). For the above $G$ and $\chi$, the generalized matrix function with respect to $G$ and $\chi$  defined on the  space of $m\times m$ matrices is
$$d^G_{\chi}(X)= \sum_{\sigma \in G}\chi(\sigma) \prod_{t=1}^m x_{t\sigma (t)}, \quad
\mbox{where $X=(x_{ij})$}. $$

If $G=S_m$
 and $\chi$ is the signum function with values $\pm 1$, then the generalized matrix function becomes the usual matrix determinant ($\det$);
   setting $\chi(\sigma)=1$ for each $\sigma \in G= S_{m}$ defines the permanent ($\per$) of the matrix.

Let $\dim \H=m$ and let $\{e_1, e_2, \dots, e_m\}$ be an orthonormal basis of $\H$. Let $P{'}$ (the transpose of $P$) be a matrix representation of a linear operator $T$ on $\H$ with respect to the basis  $\{e_1, e_2, \dots, e_m\}$. Then
(see, e.g., \cite[p.~227]{MerMul97})
$$d_{\chi}^G(P)=\frac{o(G)}{o(\chi)} \langle K(T)e^{*}, e^{*}\rangle,$$
 where $o(G)$ is the order of $G$,  $o(\chi)$ is the degree of  $\chi$,   and
  $e^{*}=e_1*e_2*\cdots *e_m$ is the decomposable symmetrized tensor of $e_1, e_2, \dots, e_m$
 (\cite[p.~155]{MerMul97}).

The following result for generalized matrix functions are immediate.

\begin{thm}\label{Thm3}
 Let $A_1, A_2, A_3$ be $m\times m$ positive semidefinite matrices. Let $G$ be a subgroup of $S_m$ and $\chi$ be an irreducible character of $G$. Then
\begin{eqnarray*}\label{dg1}
 d_{\chi}^G(A_{1}+A_{2}+A_{3})+d_{\chi}^G(A_{1})+d_{\chi}^G(A_{2})+d_{\chi}^G(A_{3})-
 \end{eqnarray*}
\vspace{-.25in}
\begin{eqnarray}\label{dg2}
\big (d_{\chi}^G(A_{1}+A_{2})+d_{\chi}^G(A_{1}+A_{3})+d_{\chi}^G(A_{2}+A_{3})\big ) \geq 0.
 \end{eqnarray}
 \end{thm}

The determinant and permanent inequalities follow at once.

\begin{cor}
Let $A_1, A_2, A_3$ be $m\times m$ positive semidefinite matrices. Then
\begin{eqnarray*}\label{det1}
 \det(A_{1}+A_{2}+A_{3})+\det(A_{1})+\det(A_{2})+\det(A_{3})-
 \end{eqnarray*}
\vspace{-.25in}
\begin{eqnarray}\label{det2}
 \big (\det(A_{1}+A_{2})+\det(A_{1}+A_{3})+\det(A_{2}+A_{3})\big ) \geq 0.
 \end{eqnarray}
\end{cor}

\begin{cor}
Let $A_1, A_2, A_3$ be $m\times m$ positive semidefinite matrices. Then
 \begin{eqnarray*}\label{per1}
 \per (A_{1}+A_{2}+A_{3})+\per(A_{1})+\per(A_{2})+\per(A_{3})-
 \end{eqnarray*}
\vspace{-.25in}
\begin{eqnarray*}\label{per2}
  \big (\per(A_{1}+A_{2})+\per(A_{1}+A_{3})+\per(A_{2}+A_{3})\big ) \geq 0.
 \end{eqnarray*}
\end{cor}

Inequality (\ref{det1}) is obtained by Lin in \cite[Theorem 1.1]{Lin20140604} by using a majorization approach,
while inequality (\ref{dg1}) in Theorem
\ref{Thm3}   confirms the strong superadditivity of the generalized matrix functions - a question raised by Lin in  \cite[Section 3]{Lin20140604}. Additionally,  our result
gives a different proof for the inequality in \cite{PTZEla14} in which an embedding approach is employed. These   results are generalizations
of the classic  inequalities $\det (A+B)\ge \det A +\det B$ and
$\per (A+B)\ge \per A +\per B$ for positive semidefinite matrices
$A$ and $B$ of the same size (see, e.g., \cite[p.~121]{BYWangMul}).

\bigskip

{\bf Acknowledgement.} The work of Haixia Chang was done during the academic year 2014-2015 when she was a Visiting Professor at Nova Southeastern University;
her work was partially supported by Scientific Research of Foundation of Shanghai Finance University (SHFUKT13-08). Fuzhen Zhang is thankful to American Institute of Mathematics (AIM) for the financial support of his attending the AIM workshop on Positivity, Oct. 13-17, 2014.
We are all grateful to Minghua Lin for sending us an early version of
\cite{Lin20140604} which initiated our work.


\end{document}